\documentclass{article}

\usepackage[repositionpullbacks,midshaft,PostScript=dvips,nohug,scriptlabels,size=2em]{diagrams}

\usepackage{latexsym}
\usepackage{amssymb}

\newcommand{\dashbk}{-}

\newcommand{\ucontents}[2]{\addcontentsline{toc}{#1}{\numberline{}{#2}}}

\newcommand{\mb}[1]{\mathbf{#1}}

\newcommand{\mr}[1]{\mathrm{#1}}

\newcommand{\fcat}[1]{\mb{#1}}

\newcommand{\slsh}{/\linebreak[0]}

\newcommand{\such}{\:|\:}
\newcommand{\without}{\setminus}

\newcommand{\implies}{\,\Rightarrow\,}

\newcommand{\Sym}{\fcat{Sym}}

\newcommand{\integers}{\mathbb{Z}}

\newarrow{Equals}=====
\newarrow{Get}....>
\newarrow{Goesto}|---{->}
\newarrow{Incl}C--->
\newarrow{Mod}--+->
\newarrow{Monic}>--->
\newarrow{NT}===={=>}

\def\today{\number\day\space \ifcase\month\or
  January\or February\or March\or April\or May\or June\or
  July\or August\or September\or October\or November\or December\fi
  \space\number\year}

\newcommand{\lwr}[1]{\mathbf{#1}}       % Object of \scat{D}
\newcommand{\complexes}{\mathbb{C}}
\newcommand{\demph}[1]{\textbf{\textup{#1}}}
\newcommand{\done}{\hfill\ensuremath{\Box}}
\newenvironment{prooflike}[1]{\begin{trivlist}\item\textbf{#1}\ }
{\end{trivlist}}
\newenvironment{proof}{\begin{prooflike}{Proof}}{\end{prooflike}}
\newcommand{\go}{\rTo\linebreak[0]}
\newcommand{\goby}[1]{\rTo^{#1}\linebreak[0]}
\newcommand{\iso}{\cong}
\newcommand{\nat}{\mathbb{N}}   
\newcommand{\url}[1]{#1}
\newcommand{\of}{\,\raisebox{0.08ex}{\ensuremath{\scriptstyle\circ}}\,}
\newcommand{\sub}{\subseteq}
\newcommand{\rationals}{\mathbb{Q}}
\newcommand{\bl}{{\scriptscriptstyle\bullet}}
\newcommand{\adj}{\mr{adj}}
\newcommand{\rmv}[2]{#1^{[#2]}}
\newcommand{\dhi}{\chi_\Sigma}
\newcommand{\pow}[1]{{[\![ #1 ]\!]}}
\newcommand{\Lau}[1]{(\!( #1 )\!)}
\newcommand{\Mat}[2]{\textrm{Mat}_{#1}(#2)}
\newcommand{\leftmat}{\left(\!\!}
\newcommand{\rightmat}{\!\!\right)}
\newcommand{\scat}[1]{\mathbf{#1}}
\newcommand{\ord}{\mathrm{o}}
\newcommand{\sgn}{\mathrm{sgn}}
\newcommand{\Fix}{\mathrm{Fix}}

\newtheorem{thm}{Theorem}[section]
\newtheorem{propn}[thm]{Proposition}
\newtheorem{lemma}[thm]{Lemma}
\newtheorem{cor}[thm]{Corollary}
\newtheorem{lotsofremarks}[thm]{Remarks}

\newtheorem{predefn}[thm]{Definition}
\newenvironment{defn}{\begin{predefn}\upshape}{\end{predefn}}
\newtheorem{preexample}[thm]{Example}
\newenvironment{example}{\begin{preexample}\upshape}{\end{preexample}}
\newenvironment{example*}[1]{\begin{preexample}\upshape}{\end{preexample}}
\newtheorem{prewarning}[thm]{Warning}

\newtheorem{preexamples}[thm]{Examples}

\title{The Euler characteristic of a category\\
as the sum of a divergent series}
\author{Tom Leinster%
\thanks{Department of Mathematics, University of Glasgow, UK;
T.Leinster@maths.gla.ac.uk.  Supported by an EPSRC Advanced Research
Fellowship.}}
\date{}

\begin{document}

\sloppy

\maketitle

\begin{abstract}
The Euler characteristic of a cell complex is often thought of as the
alternating sum of the number of cells of each dimension.  When the complex is
infinite, the sum diverges.  Nevertheless, it can sometimes be evaluated;
in particular, this is possible when the complex is the nerve of a finite
category.  This provides an alternative definition of the Euler characteristic
of a category, which is in many cases equivalent to the original one.
% ~\cite{ECC}.
\end{abstract}

% \tableofcontents

\section{Introduction}
\label{sec:intro}

What is the Euler characteristic of an infinite cell complex?  

The Euler characteristic of a finite complex is most often described as the
alternating sum of the number of cells of each dimension.  There seems little
hope of extending this formula to complexes containing infinitely many
cells of the same dimension.  However, there are interesting complexes in
which there are only finitely many cells of each dimension, but infinitely
many in total.  (The classifying space of a finite group provides an example;
see below.)  Writing $c_n$ for the number of $n$-dimensional cells, we would
like to find a sensible way of evaluating the divergent series $\sum_{n \in
\nat} (-1)^n c_n$, which could then be interpreted as the Euler characteristic
of the complex.

To see how this might work, consider a finite group $G$.  Its classifying
space $BG$ is the geometric realization of a simplicial set in which an
$n$-simplex is an $n$-tuple of elements of $G$.  The nondegenerate
$n$-simplices are the $n$-tuples of non-identity elements, so, writing
$\ord(G)$ for the order of $G$, there are $(\ord(G) - 1)^n$ of them.  A
simplicial set may be regarded as a kind of complex in which the cells
are the nondegenerate simplices, so our task is to evaluate
\begin{equation}        \label{eq:geometric}
\sum_{n \in \nat} (-1)^n (\ord(G) - 1)^n.
\end{equation}
In the Eulerian spirit of formal calculation, we apply the formula for the sum
of a geometric series, which gives the answer $1/\ord(G)$.  And indeed, it has
been established that this is the `right' value for the Euler characteristic
of $G$ (or $BG$) from several points of view: see~\cite{Wall} and~\cite{BD},
for instance.

Here is a first step towards making this rigorous.  Take a cell complex $X$
(in some non-specific sense; the exact meaning is not important for this
discussion).  Suppose that $X$ has only a finite number $c_n$ of cells of each
dimension $n$, and write $f_X(t) = \sum_{n\in\nat} c_n t^n$ for the resulting
formal power series.  It may be that $f_X$ converges in some neighbourhood of
$0$.  If so, it may also be that $f_X$ can be analytically continued to $-1$,
and it may even be that all such analytic continuations take the same value at
$-1$.  We could then define the Euler characteristic of $X$ to be that value.
Of course, if $X$ has only finitely many cells then the situation is very
simple: $f_X$ is a polynomial, there is a unique analytic continuation to
$-1$, and its value there (namely, $f_X(-1)$) is the Euler characteristic of
$X$ in the usual sense.

The purpose of this paper is to use this approach to derive a notion of the
Euler characteristic of a finite category.  (The idea was suggested to me by
Clemens Berger, to whom I am very grateful.)  This is achieved with the aid of
the nerve construction (see~\S\ref{sec:defn}), which turns a category into a
simplicial set.

We will see that when $X$ is the nerve of a finite category, the power series
$f_X$ is in fact the germ at $0$ of a rational function.  The question of
analytic continuation is then straightforward.  We can therefore carry out the
plan above to arrive at a definition of the Euler characteristic of a finite
category (at least, when the rational function does not have a pole at $-1$).
This is called the `series Euler characteristic' of the
category~(\S\ref{sec:defn}).

The new notion of series Euler characteristic is to be compared with the
original notion of the Euler characteristic of a finite category, introduced
in~\cite{ECC}.  We will see that in a large and important class of finite
categories, the two notions agree~(\S\ref{sec:positive}).  However, outside
this class the relationship breaks down: there are examples of finite
categories for which the Euler characteristic is defined in one sense but not
the other (both ways round), and further examples where both are defined
but their values differ (\S\ref{sec:negative}).

\subsubsection*{Background: the Euler characteristic of a category} 

Here are the necessary definitions from~\cite{ECC}, with notation changed
slightly to allow a matrix-based approach.

Given a natural number $m$ and a ring $K$, write $\Mat{m}{K}$ for the ring of
$m\times m$ matrices over $K$.  Given a finite category $\scat{A}$ with
objects $a_1, \ldots, a_m$, write $Z_\scat{A} \in \Mat{m}{\rationals}$ for the
matrix whose $(i, j)$-entry is the number of arrows from $a_i$ to $a_j$.  (Of
course, $Z_\scat{A}$ also depends on the order in which the objects are
listed.)  A \demph{weighting} on $\scat{A}$ is an $m$-tuple $w^\bl = (w^1,
\ldots, w^m) \in \rationals^m$ such that
\[
Z_\scat{A}
\leftmat
\begin{array}{c}
w^1     \\
\vdots  \\
w^m     
\end{array}
\rightmat
=
\leftmat
\begin{array}{c}
1       \\
\vdots  \\
1     
\end{array}
\rightmat.
\]
A \demph{coweighting} on $\scat{A}$ is an $m$-tuple $w_\bl = (w_1, \ldots,
w_m) \in \rationals^m$ such that 
\[
\leftmat 
\begin{array}{ccc}
w_1     &\cdots &w_m
\end{array}
\rightmat
Z_\scat{A}
=
\leftmat 
\begin{array}{ccc}
1       &\cdots &1
\end{array}
\rightmat.
\]
It is easy to show that if $w^\bl$ is a weighting and $w_\bl$ a coweighting on
$\scat{A}$ then $\sum_i w^i = \sum_i w_i$.  A finite category $\scat{A}$
\demph{has Euler characteristic} if it admits both a weighting and a
coweighting, and in that case its \demph{Euler characteristic}
$\chi(\scat{A})$ is $\sum_i w^i = \sum_i w_i \in \rationals$, for any
weighting $w^\bl$ and coweighting $w_\bl$.

The Euler characteristic of a category is independent of the choice of
ordering of the objects.  It is also independent of the composition and
identities; that is, it depends only on the underlying directed graph.  
(But it is usually not equal to the Euler characteristic of the underlying
graph, `vertices minus edges'.  The Euler characteristics of categories and
graphs are compatible in a different sense: Proposition~2.10
of~\cite{ECC}.)

An important special case is when $Z_\scat{A}$ is invertible.  Then $\scat{A}$
is said to have \demph{M\"obius inversion}, there are a unique weighting and a
unique coweighting, and $\chi(\scat{A})$ is the sum of the entries of
$Z_\scat{A}^{-1}$. 

The Euler characteristic of categories enjoys many good properties.  It is
invariant under equivalence and behaves predictably with respect to products,
fibrations, etc.  It is also compatible with Euler characteristics of other
types of structure, including topological spaces, graphs, posets, groups,
manifolds, and orbifolds.

\subsubsection*{Related work}

A different notion of M\"obius inversion for categories is considered
in~\cite{CLL} and~\cite{Hai}.  See \S 4 of~\cite{ECC} for a discussion of the
relationship.

The observation that the Euler characteristic of a finite group can be
calculated by formal summation of the geometric series~(\ref{eq:geometric}) is
probably nearly as old as the concept of the classifying space of a group.  (I
learned it from a talk of John Baez.)  A group can be regarded as a one-object
category in which all morphisms are invertible; from that viewpoint, our
purpose is to extend this formal method from finite groups to finite
categories.

Further references to related work can be found in~\cite{ECC}.

\section{The series Euler characteristic of a category}
\label{sec:defn}

In this section we define the series Euler characteristic of a finite
category, see how the definition works in the motivating case where the
category is a group, and find a way to calculate it.

First let us recall some facts about formal power series.  For any field $K$,
there is a commutative diagram
\[
\begin{diagram}
K[t]    &\rIncl         &K\pow{t}       \\
\dIncl  &               &\dIncl         \\
K(t)    &\rIncl         &K\Lau{t}         \\
\end{diagram}
\]
of rings.  Here $K[t]$ is the ring of polynomials over $K$ and $K\pow{t}$ is
the ring of formal power series $\sum_{n \in \nat} a_n t^n$ ($a_n \in K$).
These are both integral domains, and their fields of fractions are in the
bottom row: the field $K(t)$ of rational expressions over $K$, and the field
$K\Lau{t}$ of formal Laurent series over $K$ (expressions
$\sum_{n\in\integers} a_n t^n$ such that $\{n \leq 0 \such a_n \neq 0\}$ is
finite).

The canonical inclusions of $K\pow{t}$ and $K(t)$ into $K\Lau{t}$ make it
possible to ask whether a formal power series `is rational'; in other words,
whether the element of $K\Lau{t}$ corresponding to the power series is in the
image of $K(t)$.  When $K$ is a subfield of $\complexes$, the following
analytic criterion applies.  Let $f \in K\pow{t}$ and let $p, q \in K[t]$ be
coprime polynomials.  Then $f = p/q$ in $K\Lau{t}$ if and only if there is a
neighbourhood $U$ of $0$ in $\complexes$ such that for all $z \in U$, $q(z)
\neq 0$ and $f(z)$ converges to $p(z)/q(z)$.

We will also need some notation for matrices.  Let $m \in \nat$ and let $K$ be
a commutative ring.  We write $s: \Mat{m}{K} \go K$ for the $K$-linear map
defined by $s(M) = \sum_{i, j} M_{ij}$.  Every matrix $M \in \Mat{m}{K}$ has
an \demph{adjugate} $\adj(M) \in \Mat{m}{K}$, defined by
\[
(\adj(M))_{ij} 
=
(-1)^{i + j}
\cdot
\det( M \textrm{ with its } j \textrm{th row and } i \textrm{th column
deleted} )
\]
and satisfying
\[
M \cdot \adj(M) = \adj(M) \cdot M = \det(M) \cdot I.
\]

\begin{lemma}
Let $M$ be a square matrix over a field $K$.  Then $\sum_{n\in\nat} s(M^n) t^n
\in K\pow{t}$ is rational.
\end{lemma}

\begin{proof}
Write
\[
F(t)
=
\sum_{n\in\nat} M^n t^n 
\in 
\Mat{m}{K\pow{t}}.
\]
Then
\[
(I - Mt) F(t) 
=
I,
\]
so 
\[
\det(I - Mt) \cdot F(t)
=
\adj(I - Mt).
\]
Applying the $K\pow{t}$-linear map $s: \Mat{m}{K\pow{t}} \go K\pow{t}$,
\[
\det(I - Mt) \cdot s(F(t))
=
s(\adj(I - Mt)).
\]
But $s(F(t)) = \sum s(M^n) t^n$, and $\det(I - Mt)$ is not the zero polynomial
(since its value at $t = 0$ is $1$), so $\sum s(M^n) t^n$ is rational and
equal to
\[
\frac{s(\adj(I - Mt))}{\det(I - Mt)}
\in
K(t).
\]
\done
\end{proof}

Given a simplicial set $X$ with only finitely many simplices of each
dimension, let $c_n$ be the number of nondegenerate $n$-simplices and $f_X(t)
= \sum_{n\in\nat} c_n t^n \in \rationals\pow{t}$.  Recall that the
\demph{nerve} $N\scat{A}$ of a category $\scat{A}$ is a simplicial set in
which an $n$-simplex is a chain
\begin{equation}        \label{eq:nerve}
x_0 \goby{\phi_1} x_1 \goby{\phi_2} \cdots \goby{\phi_n} x_n
\end{equation}
of arrows in $\scat{A}$; such an $n$-simplex is degenerate if and only if some
$\phi_i$ is an identity.  When $\scat{A}$ is finite, write $f_\scat{A} =
f_{N\scat{A}} \in \rationals\pow{t}$.

\begin{thm}
For any finite category $\scat{A}$, the formal power series $f_\scat{A}$ is
rational (over $\rationals$).
\end{thm}

\begin{proof}
Order the objects of $\scat{A}$ as $a_1, \ldots, a_m$ and let $Z_\scat{A}$ be
the matrix of $\scat{A}$ with respect to this ordering, as
in~\S\ref{sec:intro}.  For each $i$ and $j$, the number of non-identity arrows
from $a_i$ to $a_j$ is $(Z_\scat{A} - I)_{ij}$.  The number of nondegenerate
$n$-simplices~(\ref{eq:nerve}) beginning at $a_i$ and ending at $a_j$ is,
therefore, $((Z_\scat{A} - I)^n)_{ij}$.  Hence the total number $c_n$ of
nondegenerate $n$-simplices is $s((Z_\scat{A} - I)^n)$.  The result now
follows from the Lemma.  \done
\end{proof}

The series $\sum (-1)^n c_n$ is in general divergent.
(Proposition~2.11 of~\cite{ECC} gives exact conditions for it to
converge.)  Nevertheless, the Theorem provides a way to `evaluate' it,
returning an answer $f_\scat{A}(-1) \in \rationals \cup \{\infty\}$.

\begin{defn}
A finite category $\scat{A}$ \demph{has series Euler characteristic} if
$f_\scat{A}(-1) \in \rationals$.  In that case, its \demph{series Euler
characteristic} is $\dhi(\scat{A}) = f_\scat{A}(-1)$.
\end{defn}

The proofs tell us that
\[
f_\scat{A}(t)
=
\frac{s(\adj(I - (Z_\scat{A} - I)t))}{\det(I - (Z_\scat{A} - I)t)}.
\]

\begin{example*}{Monoids}
Let $G$ be a monoid of finite order $\ord(G)$, and let $\scat{A}$ be the
corresponding one-object category.  Then $\sum (-1)^n c_n = \sum (-1)^n
(\ord(G) - 1)^n$ is divergent (unless $G$ is trivial), but the rational
function
\[
f_\scat{A}(t)
=
\sum_{n\in\nat} c_n t^n
=
\frac{1}{1 - (\ord(G) - 1)t}
\]
has value $1/\ord(G)$ at $t = -1$.  Hence $\dhi(\scat{A}) = 1/\ord(G)$.
\end{example*}

A change of variable will be useful.  Put $u = 1 + 1/t$ and write
\[
g_\scat{A}(u) 
=
\frac{s(\adj(Z_\scat{A} - uI))}{\det(Z_\scat{A} - uI)}
\in \rationals(u).
\]
Then $f_\scat{A}(t) = (1 - u) g_\scat{A}(u)$.  Hence $\dhi(\scat{A}) =
g_\scat{A}(0)$, with one side defined if and only if the other is.

We will need to be able to compute values of $\dhi$.  The key observation is
that for an $m \times m$ matrix $M$ over a commutative ring,
% 
% \begin{equation}        \label{eq:det-like}
\[
s(\adj(M)) 
=
\sum_{\sigma \in S_m} 
\sgn(\sigma) \cdot
F((M_{i, \sigma(i)})_{i \in \lwr{m}})
\]
% \end{equation}
% 
where $\lwr{m} = \{1, \ldots, m\}$ and $F$ is the symmetric function defined
by 
\[
F((x_i)_{i \in I})
=
\sum_{i \in I}
\prod_{j \in I\without\{i\}} 
x_j
\]
for any finite family $(x_i)_{i \in I}$.  This follows from the analogous
formula for determinants.  

Given an $m \times m$ matrix $M$ and a subset $R$ of $\lwr{m}$, denote by
$\rmv{M}{R}$ the matrix obtained from $M$ by deleting the $i$th row and
$i$th column for every $i \in R$.  Write $\Sym(S)$ for the group of
permutations of a set $S$. 

\begin{propn}   \label{propn:comp}
Let $m \in \nat$ and let $Z$ be an $m \times m$ matrix over a commutative
ring.  Then
\begin{equation}        \label{eq:det}
\det(Z - uI)
=
\sum_{r = 0}^m
(-1)^r
d_r
u^r
\quad
\textrm{where}
\quad
d_r
=
\sum_{R \sub \lwr{m},\  |R| = r}
\det(\rmv{Z}{R})
\end{equation}
and 
\begin{equation}        \label{eq:sadj}
s(\adj(Z - uI))
=
\sum_{r = 0}^m
(-1)^r
e_r 
u^r
\quad
\textrm{ where }
\quad
e_r
=
\sum_{R \sub \lwr{m},\  |R| = r}
s(\adj(\rmv{Z}{R})).
\end{equation}
\end{propn}

\begin{proof}
Equation~(\ref{eq:det}) is classical~\cite{Jac}.  For~(\ref{eq:sadj}), 
if $\sigma \in S_m$ then, by a short calculation,
% 
% \begin{equation}        \label{eq:sym-sum}
\[
F(((Z - uI)_{i, \sigma(i)})_{i \in \lwr{m}})
=
\sum_{R \sub \Fix(\sigma)}
(-u)^{|R|}
F((Z_{i, \sigma(i)})_{i \in \lwr{m}\without R}).
\]
% \end{equation}
% 
Hence 
\begin{eqnarray*}
s(\adj(Z - uI)) &=      &
\sum_{\sigma \in S_m}
\sgn(\sigma) \cdot
F(((Z - uI)_{i, \sigma(i)})_{i \in \lwr{m}})    
\label{eq:comp-1}       \\
        &=      &
\sum_{\sigma \in S_m, \ R \sub \Fix(\sigma)}
\sgn(\sigma) \cdot
(-u)^{|R|} 
F((Z_{i, \sigma(i)})_{i \in \lwr{m}\without R})
\label{eq:comp-2}       \\
        &=      &
\sum_{R \sub \lwr{m}} 
(-u)^{|R|}
\sum_{\sigma' \in \Sym(\lwr{m}\without R)} 
\sgn(\sigma') \cdot
F((Z_{i, \sigma'(i)})_{i \in \lwr{m}\without R})
\label{eq:comp-3}       \\
        &=      &
\sum_{R \sub \lwr{m}}
(-u)^{|R|} s(\adj(\rmv{Z}{R})),
\label{eq:comp-4}
\end{eqnarray*}
as required.  (In fact, the same argument proves~(\ref{eq:det}) too, by
changing $s(\adj(\dashbk))$ to $\det$ and $F$ to product throughout.)
\done
\end{proof}

Given a finite category $\scat{A}$, take $Z = Z_\scat{A}$ and write
$d^\scat{A}_r = d_r$, $e^\scat{A}_r = e_r$.  Denote by $l$ the least number
such that $d^\scat{A}_l \neq 0$.  Then $\scat{A}$ has series Euler
characteristic if and only if $e^\scat{A}_r = 0$ for all $r < l$, and in that
case,
\begin{equation}        \label{eq:dhi-comp}
\dhi(\scat{A}) 
=
e^\scat{A}_l/d^\scat{A}_l. 
\end{equation}

\section{Positive results}
\label{sec:positive}

Among the finite categories that have Euler characteristic, those
with M\"obius inversion form an important class.  Any finite
category equivalent to one with M\"obius inversion also has Euler
characteristic, and this larger class encompasses most of the
examples in~\cite{ECC}: finite monoids, groupoids, posets,
categories in which the only endomorphisms are automorphisms (or
equivalently, the only idempotents are identities), and
categories admitting an epi-mono factorization system.  We show
that Euler characteristic and series Euler characteristic agree
on this class.

\begin{lemma}   \label{lemma:repetition}
Let $M$ be a square matrix such that for some $i \neq
j$, the $i$th and $j$th rows are equal and the $i$th and $j$th columns are
equal.  Then $s(\adj(M)) = 0$.
\end{lemma}

\begin{proof} 
First suppose that $i = 1$ and $j = 2$.  Then every entry of $\adj(M)$ is zero
except perhaps for the four in the top-left corner, which are
\[
\leftmat
\begin{array}{cc}
y       &-y     \\
-y      &y      
\end{array}
\rightmat
\]
where $y$ is the $(1, 1)$-minor of $M$.  Hence $s(\adj(M)) = 0$.

The general case is handled similarly.  Alternatively, it may be reduced to
the case $(i, j) = (1, 2)$ by showing that $s(\adj(M))$ is unchanged when a
permutation is applied to both the rows and the columns of $M$.
\done
\end{proof}

\begin{thm} 
Let $\scat{A}$ be a finite category equivalent to some category with M\"obius
inversion.  Then $\scat{A}$ has both Euler characteristic and series Euler
characteristic, and $\chi(\scat{A}) = \dhi(\scat{A})$.
\end{thm}

\begin{proof} 
Order the objects of $\scat{A}$ so that the isomorphism
classes are grouped together: 
\[
a_1^1, \ldots, a_1^{q_1}, \ldots, 
a_n^1, \ldots, a_n^{q_n}
\]
where $a_i^j \iso a_{i'}^{j'}$ if and only if $i = i'$, and where
each $q_i$ is nonzero.  Let $\scat{B}$ be the full subcategory
$\{ a_1^1, a_2^1, \ldots, a_n^1 \}$, a skeleton of $\scat{A}$.
Now $\scat{A}$ is equivalent to some category $\scat{B}'$ with
M\"obius inversion, and any category with M\"obius inversion is
skeletal, so $\scat{B}'$ is isomorphic to $\scat{B}$ and $\scat{B}$ has
M\"obius inversion.
% 
%  and is
% therefore (up to isomorphism) the unique skeletal category equivalent to
% $\scat{A}$.  Categories with M\"obius inversion are skeletal, so $\scat{B}$
% has M\"obius inversion.  
Hence $\chi(\scat{B})$ is defined.  Since Euler
characteristic is invariant under equivalence, $\chi(\scat{A})$ is also
defined and $\chi(\scat{A}) = \chi(\scat{B})$.

Let $R \sub \{1, \ldots, q_1 + \cdots + q_n\}$.  By
Lemma~\ref{lemma:repetition}, $s(\adj(\rmv{Z_\scat{A}}{R})) = 0$ unless $R$
omits at most one element of each isomorphism class, and in particular has at
least $l = q_1 + \cdots + q_n - n$ elements.  So $e^\scat{A}_r = 0$ for all $r
< l$.  If $R$ has $l$ elements then in order for
$s(\adj(\rmv{Z_{\scat{A}}}{R}))$ to be nonzero, $R$ must omit exactly one
element of each isomorphism class, in which case $\rmv{Z_\scat{A}}{R} =
Z_\scat{B}$.  Hence
\[
e^{\scat{A}}_l = q_1 \cdots q_n s(\adj(Z_\scat{B})).
\]
Similarly, $d^\scat{A}_r = 0$ for all $r < l$ and 
\[
d^{\scat{A}}_l = q_1 \cdots q_n \det(Z_\scat{B}).
\]
But $\scat{B}$ has M\"obius inversion, that is, $\det(Z_{\scat{B}}) \neq
0$, so $d^\scat{A}_l \neq 0$.  Hence 
\[
\dhi(\scat{A})
=
e^\scat{A}_l / d^\scat{A}_l
=
s(\adj(Z_\scat{B})) / \det(Z_\scat{B})
=
s(Z_\scat{B}^{-1})
=
\chi(\scat{B})
=
\chi(\scat{A}),
\]
using~(\ref{eq:dhi-comp}) in the first equality.
\done
\end{proof}

The next result gives a class of categories that need not have Euler
characteristic, but do have series Euler characteristic.

\begin{propn}   \label{propn:diag} 
If $\scat{A}$ has either a weighting or a coweighting and
$Z_\scat{A}$ is diagonalizable then $\scat{A}$ has series Euler
characteristic. 
\end{propn}

\begin{proof} 
We may write $Z_\scat{A} = PDP^{-1}$, with $D$ the diagonal matrix on 
$(\lambda_1, \ldots, \lambda_m)$, and by duality we may assume that $\scat{A}$
has a coweighting $w_\bl$.  

For $n \in \nat$, 
\[
s((Z_\scat{A} - I)^n)
=
\sum_{i, j, k \in \lwr{m}} 
P_{ij} (\lambda_j - 1)^n (P^{-1})_{jk}
=
\sum_{j \in \lwr{m}} p_j p'_j (\lambda_j - 1)^n
\]
where $p_j$ is the $j$th column-sum of $P$ and $p'_j$ the $j$th row-sum of
$P^{-1}$.  Hence
\begin{equation}        \label{eq:power-series}
f_\scat{A}(t)
=
\sum_{j \in \lwr{m}} \frac{p_j p'_j}{1 - (\lambda_j - 1)t}.
\end{equation}
It suffices to prove that $p_j p'_j = 0$ for all $j$ such that
$\lambda_j = 0$.  Indeed, suppose that $\lambda_j = 0$.  Then,
writing $P_j$ for the $j$th column of $P$, we have
\[
p_j 
=
\leftmat
\begin{array}{ccc}
1       &\cdots         &1
\end{array}
\rightmat
P_j
=
w_\bl Z_\scat{A} P_j
=
w_\bl \lambda_j P_j
=
0,
\]
as required.
\done
\end{proof}

\begin{example} \label{eg:four-objs} 
Let $\scat{A}$ be the $4$-object category in Example~1.11(d) of~\cite{ECC},
which admits a coweighting (since it has an initial object) but no weighting,
and so does not have Euler characteristic.  Then
\[
Z_\scat{A} 
=
\leftmat
\begin{array}{cccc}
2       &2      &1      &1      \\
2       &2      &1      &2      \\
1       &1      &1      &1      \\
0       &0      &0      &1      
\end{array}
\rightmat,
\]
which is diagonalizable, so $\scat{A}$ has series Euler characteristic.  (In
fact, it can be shown using~(\ref{eq:power-series}) that $\dhi(\scat{A}) =
1$.)
\end{example}

\section{Negative results}
\label{sec:negative}

We have already seen that Euler characteristic and series Euler
characteristic are defined and agree in a large and important
class of finite categories, namely, those equivalent to some
category with M\"obius inversion.  In this section we see that
outside this class, the relationship breaks down.

We first see that the properties of having Euler characteristic and having
series Euler characteristic are logically independent.  In other words, all
four possibilities occur: a category may have both Euler characteristic and
series Euler characteristic (as in \S\ref{sec:positive}), series Euler
characteristic but not Euler characteristic (Example~\ref{eg:four-objs}),
Euler characteristic but not series Euler characteristic
(Example~\ref{eg:old-not-new}), or neither (Example~\ref{eg:neither}).
Furthermore, even when both are defined, they do not necessarily agree
(Example~\ref{eg:incompatible}).

Given the disagreement between the two definitions, one might ask which (if
either) should be regarded as superior.  A major defect of series Euler
characteristic is that it is not invariant under equivalence
(Example~\ref{eg:not-invariant}).  In contrast, ordinary Euler characteristic
is invariant not only under equivalence but also under the existence of an
adjunction (Proposition~2.4 of~\cite{ECC}).  I do not know whether series
Euler characteristic enjoys the same properties with respect to products,
fibrations, etc.

For the examples, we will need to know something about which matrices arise
from categories.  Let us say that a square matrix $Z$ of natural numbers is
\demph{the matrix of a category} if there exists a finite category $\scat{A}$
such that $Z = Z_\scat{A}$ (with respect to some ordering of the objects).
Any such matrix $Z$ is certainly \demph{reflexive} ($Z_{ii} \geq 1$ for all
$i$) and \demph{transitive} ($Z_{ij}, Z_{jk} \geq 1 \implies Z_{ik} \geq 1$).
These necessary conditions are not sufficient; for instance, it can be shown
that
\[
\leftmat
\begin{array}{cc}
1       &2      \\
1       &2      
\end{array}
\rightmat
\]
is not the matrix of a category.  However, we do have:

\begin{lemma} 
Let $Z$ be a transitive square matrix of natural numbers whose diagonal
entries are all at least $2$.  Then $Z$ is the matrix of a category.
\end{lemma}

\begin{proof} 
Suppose that $Z$ is an $m\times m$ matrix.  We define a category structure on
the directed graph with objects $1, \ldots, m$ and with $Z_{ij}$ arrows from
$i$ to $j$, for each $i$ and $j$.  For each $i$, choose an arrow $1_i: i \go
i$.  For each pair $(i, j)$ such that $Z_{ij} \geq 1$, choose an arrow
$\phi_{ij}: i \go j$, with $\phi_{ii} \neq 1_i$ for all $i$.  To define
composition, take arrows $i \goby{\alpha} j \goby{\beta} k$.  If either
$\alpha$ or $\beta$ is an identity, it is clear how $\beta \of \alpha$ must be
defined; otherwise, put $\beta \of \alpha = \phi_{ik}$.  \done
\end{proof}

\begin{cor} 
Let $Z$ be a square matrix of positive integers whose diagonal entries are all
at least $2$.  Then $Z$ is the matrix of a category.  \done
\end{cor}
 
All of the examples that follow use this corollary without
mention.  They can be verified using Proposition~\ref{propn:comp}
and the remark after it.

\begin{example*}{New doesn't contain old}       \label{eg:old-not-new} 
A category may have Euler characteristic but not series Euler characteristic.
For example, there is a category $\scat{A}$ with
\[
Z_\scat{A}
=
\leftmat
\begin{array}{cccc}
6       &6      &15     &9      \\
6       &6      &6      &6      \\
6       &6      &9      &7      \\
6       &30     &9      &15
\end{array}
\rightmat,
\qquad
g_\scat{A}(u) 
=
\frac{4(1 + u)}{u(36 - u)}.
\]
Then $\scat{A}$ has a weighting $(1/6, 0, 0, 0)$ and a coweighting $(0, 1/6,
0, 0)$, so $\chi(\scat{A}) = 1/6$.  But $g_\scat{A}$ has a pole at $0$, so
$\dhi(\scat{A})$ is undefined.
\end{example*}

\begin{example*}{Union isn't everything}        \label{eg:neither}
A finite category may have neither Euler characteristic nor series
Euler characteristic.  For example, there is a category $\scat{A}$
with
\[
Z_\scat{A} =
\leftmat
\begin{array}{cc}
2       &4      \\
1       &2      \\
\end{array}
\rightmat,
\qquad
g_\scat{A}(u) 
= 
\frac{1 + 2u}{u(4 - u)}.
\]
Then $\scat{A}$ does not have a weighting \emph{or} a coweighting, so
certainly does not have Euler characteristic; $\dhi(\scat{A})$ is
also undefined.  
\end{example*}

\begin{example*}{Disagreement on intersection}  \label{eg:incompatible}
Even if a category has Euler characteristic in both senses, their values may
differ.  For example, there is a category $\scat{A}$ with
\[
Z_\scat{A}
=
\leftmat
\begin{array}{ccc}
2       &2      &2      \\
2       &2      &2      \\
2       &8      &5
\end{array}
\rightmat,
\qquad
g_\scat{A}(u)
=
\frac{3}{9 - u}.
\]
Then $\chi(\scat{A}) = 1/2$ (since $(1/2, 0, 0)$ is both a weighting and a
coweighting) but $\dhi(\scat{A}) = 1/3$.
\end{example*}

\begin{example*}{Not equivalence-invariant} \label{eg:not-invariant} 
Series Euler characteristic is not invariant under equivalence.  For example,
we may choose a category $\scat{A}$ satisfying
\[
Z_\scat{A} 
=
\leftmat
\begin{array}{cc}
3       &3      \\
2       &2      
\end{array}
\rightmat,
\qquad
g_\scat{A}(u) = \frac{2}{5 - u},
\qquad
\dhi(\scat{A}) = \frac{2}{5}.
\]
Write the objects of $\scat{A}$ as $a_1, a_2$ and form a new category
$\scat{B}$ by adjoining an object $a_3$ isomorphic to $a_2$.  Then
$\scat{B}$ is equivalent to $\scat{A}$ and
\[
Z_\scat{B}
=
\leftmat
\begin{array}{ccc}
3       &3      &3      \\
2       &2      &2      \\
2       &2      &2 
\end{array}
\rightmat,
\qquad
g_\scat{B}(u) = \frac{3}{7 - u},
\qquad
\dhi(\scat{B}) = \frac{3}{7}.
\]
\end{example*}

The remaining examples concern Proposition~\ref{propn:diag}, which gives
sufficient conditions for series Euler characteristic to be defined.

\begin{example}
The Proposition is sharp, in the sense that neither of its hypotheses can be
dropped.  Example~\ref{eg:neither} shows that we cannot drop the first
hypothesis (that $\scat{A}$ admits a weighting or a coweighting), since there
$Z_\scat{A}$ is diagonalizable but $\dhi(\scat{A})$ is undefined.  To see that
we cannot drop the second hypothesis (diagonalizability of $Z_\scat{A}$),
take the following example:
\[
Z_\scat{A} 
= 
\leftmat
\begin{array}{ccc}
2       &3      &5      \\
2       &3      &5      \\
2       &1      &3      
\end{array}
\rightmat,
\qquad
g_\scat{A}(u)
=
\frac{2 + 3u}{u(8 - u)}.
\]
Then $\scat{A}$ has a weighting, but $\dhi(\scat{A})$ is undefined.
\end{example}

\begin{example}
Even when $\scat{A}$ does have a weighting and $Z_\scat{A}$ is diagonalizable,
$\dhi(\scat{A})$ need not be the total weight $\sum_i w^i$ of every weighting
$w^\bl$.  Indeed, the total weight may vary with the weighting chosen.  For
example, there is a category $\scat{A}$ with
\[
Z_\scat{A}
=
\leftmat
\begin{array}{cc}
2       &3      \\
2       &3      
\end{array}
\rightmat,
\qquad
g_\scat{A}(u) = \frac{2}{5 - u},
\qquad
\dhi(\scat{A}) = \frac{2}{5}.
\]
Then $Z_\scat{A}$ is diagonalizable, and $(1/2, 0)$ and $(1/3, 0)$ are
weightings whose total weights are different.
\end{example}

\paragraph{Acknowledgements}   I owe many thanks to Clemens Berger for the
basic idea of this paper, that the Euler characteristic of a category
$\scat{A}$ might be the value at $-1$ of some analytic continuation of the
power series $f_\scat{A}$.  He made this suggestion during a month-long visit
that I made to the University of Nice in the spring of 2007.  Much of this
work was carried out there, and I am very grateful to Eugenia Cheng for
inviting me.  I also thank Arnaud Beauville, Christina Cobbold and Andr\'e
Hirschowitz.


\small
\begin{thebibliography}{MMM}
\ucontents{section}{References}

\bibitem[BD]{BD}
John Baez, James Dolan,
From finite sets to Feynman diagrams,
in \emph{Mathematics Unlimited---2001 and Beyond},  
Springer, 2001.

\bibitem[CLL]{CLL}
Mireille Content, Fran\c{c}ois Lemay, Pierre Leroux,
Cat\'egories de M\"obius et fonctorialit\'es: un cadre
g\'en\'eral pour l'inversion de M\"obius, 
\emph{J.\ Combin.\ Theory Ser.\ A} 28 (1980), 169--190.

\bibitem[H]{Hai}
John Haigh,
On the M\"obius algebra and the Grothendieck ring of a finite category,
\emph{J.\ London Math.\ Soc.\ (2)} 21 (1980), 81--92.

\bibitem[J]{Jac}
Nathan Jacobson, 
\emph{Lectures in Abstract Algebra, Vol.~II: Linear Algebra}, 
Van Nostrand, 1953.

\bibitem[L]{ECC}
Tom Leinster,
The Euler characteristic of a category, 
\url{math.CT\slsh 0610260} (2006).

\bibitem[W]{Wall}
C.T.C. Wall,
Rational Euler characteristic,
\emph{Proc.\ Cambridge Philos.\ Soc.}\ 
57 (1961), 182--184.

\end{thebibliography}
\end{document}